\documentclass[10pt]{amsart}

\newtheorem{theorem}{Theorem}[section]
\newtheorem{lemma}[theorem]{Lemma}

\newtheorem{corollary}[theorem]{Corollary}

\usepackage{graphicx,bm}
\usepackage{tikz}
\usepackage{verbatim}
\usepackage{amssymb,amsmath,amsthm}
\usepackage{amsfonts}
\usepackage{latexsym}
\usepackage[mathscr]{eucal}
\usepackage{setspace} 
\newcommand{\lo}{\ensuremath{\lambda_1}}
\newcommand{\lt}{\ensuremath{\lambda_2}}
\newcommand{\la}{\ensuremath{\lambda}}

\begin{document}
\pagestyle{plain}

\title{Fixed Block Configuration Group Divisible Designs with Block Size 6}
\author{  M.S. Keranen and M.R. Laffin }%\date{\today}
%\begin{linenumbers}

\begin{abstract}
 
We present constructions and results about GDDs with two groups and block size
6.  We study those GDDs in which each block has configuration $(s,t)$, that
is in which each block has exactly $s$ points from one of the two groups and $t$
points from the other.  We show the necessary conditions
are sufficient for the existence of $GDD(n,2,6; \lambda_1, \lambda_2)$s with
fixed block configuration $(3,3)$.  For configuration $(1,5)$,
we give minimal or near-minimal index examples for all group sizes $n \geq 5$
except $n=10, 15, 160$, or $190$.  For configuration $(2,4)$,
we provide constructions for several families of
$GDD(n,2,6;\lambda_1,\lambda_2)$s.\footnote{Submitted for publication to
Discrete Mathematics.}
\end{abstract}

\maketitle 
\section{Introduction}

A group divisible design GDD($n,m,k$; \lo, \lt) is a collection of $k$ element
subsets of a $v$-set {\bf X} called blocks which satisfies the following
properties: each point of {\bf X} appears in the same number, $r$, of the $b$
blocks; the $v=nm$
elements of {\bf X} are partitioned into $m$ subsets (called groups) of size $n$
each; pairs of points within the same group are called first associates of each
other and appear in \lo\, blocks; pairs of points not in the same group are
second associates and appear in \lt\, blocks together. If we require that $m=2$
and each block intersects one group in
$s$ points and $t=k-s$ points in the other, we say the design has fixed  block 
configuration ($s,t$).\\ \\
%%% history
In \cite{3gdd} the authors settled the existence for group divisible designs
with block size three and first and second associates, $m$ groups of size $n$
where $m,n\geq 3$. The problem of finding necessary and sufficient conditions
for $m=2$ or $v=2n$ and block size four was established in  \cite{k4oddandeven}.
In \cite{k42grps}, the necessary conditions are shown to be sufficient for
$3\leq n \leq 8$.  New conditions and results were presented in \cite{k43grps}
with three groups and block size four, in particular, constructions were given
to show that the necessary conditions are sufficient for all GDDs with three
groups and group sizes two, three, and five with two exceptions.  In
\cite{5gdd23}, Hurd, Mishra and Sarvate gave new results for general fixed block
configuration GDD$(n,2,k;\lo,\lt)$, as well as new necessary and sufficient
conditions for $k=5$ and configuration (2,3). Hurd and Sarvate in \cite{5gdd14}
gave similar results for $k=5$ and configuration $(1,4)$. 
 Unless
otherwise stated, $m=2$ is assumed from now on. \\ \\
The purpose of this article is to establish similar results for GDDs with block
size six and two groups. In this paper, we consider each possible configuration
type: (3,3), (2,4)
and
(1,5).

\subsection{Necessary Conditions}
For GDDs with block size six and two groups there are two necessary conditions
on the number of blocks $b$, and the number of blocks a point appears in $r$. 
\begin{theorem}\label{thm1} The following conditions are necessary for the
existence of a GDD$($n, $2$,$6$;\lo,\lt$)$. 
\begin{enumerate}
\item The number of blocks is
$b=\displaystyle{\frac{\lambda_1(n)(n-1)+\lambda_2n^2}{15}}$.
\item The number of blocks a point appears in is
$r=\displaystyle{\frac{\lambda_1(n-1)+\lambda_2n}{5}}$.
\end{enumerate}
\end{theorem}

\begin{proof}
\begin{enumerate}
For condition (1), we count the total number of blocks, $b$. Each block has ${6
\choose{2}}=15$ pairs. Thus the total number of blocks must be divisible by 15.
Consider a point $v$.  There are exactly $\lambda_1(n-1)$ pairs containing
another point from the same group, and  $\lambda_2 n$ pairs with a point from
the other group. Thus the total number of pairs is
$15b=\lambda_1(n)(n-1)+\lambda_2n^2$ and the total number of blocks is
$b=\displaystyle{\frac{\lambda_1(n)(n-1)+\lambda_2n^2}{15}}$. %%% part 2
For condition (2), consider a point $v$. In any block with $v$ there are 5 pairs
containing $v$ and thus the total number of blocks containing $v$ must be
divisible by 5. Further $v$ appears in a block $\lambda_1$ times with every
other point in its same group, which is $n-1$ points, and it appears $\lambda_2$
times with every point in the other group ($n$ points in the other group). Thus
the total number of blocks that $v$ appears in is
$r=\displaystyle{\frac{\lambda_1(n-1)+\lambda_2n}{5}}$. 
\end{enumerate}
\end{proof}

\noindent These two necessary conditions on $b$ and $r$ determine possibilities
for the parameter $n$ and the indices $\lambda_1$ and $\lambda_2$. Table
\ref{possiblen} summarizes this relationship.
\begin{center}
\begin{table}
\caption{Possible values of $n$ with respect to $\lo,\lt$}
\scalebox{0.7}{
\begin{tabular}{l||r|r|r|r|r}
 \hline%%0
 \multicolumn{1}{c|}{$\pmod{15}$} &$\lambda_1\equiv 0 \mod 5$ &$\lambda_1\equiv
1\mod 5$&$\lambda_1\equiv 2\mod 5$&$\lambda_1\equiv 3\mod 5$&$\lambda_1\equiv
4\mod 5$  \\
\hline\hline%%%%%%%%%%%%%%%%%%%%%%%%%%%%%%%%%%%%%%%%%%%
 \multicolumn{1}{c|}{$\lt\equiv 0$} &Any $n$ &$n\equiv 1\mod 5$&$n\equiv 1 \mod
5$&$n\equiv 1 \mod 5$&$n\equiv 1 \mod 5$  \\
 \hline%%%%%%%%%%%%%%%%%%%%%%%%%%%%%%%%%%%%%%%%%%%
 \multicolumn{1}{c|}{$\lt\equiv 1$} &impossible  &$n\equiv 3,8\mod 15$&$n\equiv
9,14 \mod 15$&$n\equiv 2,12 \mod 15$&impossible  \\
 \hline%%%%%%%%%%%%%%%%%%%%%%%%%%%%%%%%%%%%%%%%%%%
 \multicolumn{1}{c|}{$\lt\equiv 2$} &impossible  &$n\equiv 12\mod 15$&$n\equiv
3,8 \mod 15$&impossible&$n\equiv 9 \mod 15$\\
 \hline%%%%%%%%%%%%%%%%%%%%%%%%%%%%%%%%%%%%%%%%%%%
 \multicolumn{1}{c|}{$\lt\equiv 3$} &$n\equiv 0\mod 5$  &$n\equiv 4\mod
5$&impossible&$n\equiv 3\mod 5$&$n\equiv 2 \mod 5$\\
 \hline%%%%%%%%%%%%%%%%%%%%%%%%%%%%%%%%%%%%%%%%%%%
 \multicolumn{1}{c|}{$\lt\equiv 4$} &$n\equiv 0\mod 15$  &impossible&$n\equiv
2,12\mod 15$&$n\equiv 9,14\mod 15$&$n\equiv 3,8 \mod 15$\\
 \hline%%%%%%%%%%%%%%%%%%%%%%%%%%%%%%%%%%%%%%%%%%%
 \multicolumn{1}{c|}{$\lt\equiv 5$} &$n\equiv 0\mod 3$  &$n\equiv 6 \mod
15$&$n\equiv 6,11\mod 15$&$n\equiv 6,11\mod 15$&$n\equiv 6,11 \mod 15$\\
 \hline%%%%%%%%%%%%%%%%%%%%%%%%%%%%%%%%%%%%%%%%%%%
 \multicolumn{1}{c|}{$\lt\equiv 6$} &$n\equiv 0\mod 5$  &$n\equiv 3 \mod
5$&$n\equiv 4\mod 5$&$n\equiv 2\mod 5$&impossible\\
 \hline%%%%%%%%%%%%%%%%%%%%%%%%%%%%%%%%%%%%%%%%%%%
 \multicolumn{1}{c|}{$\lt\equiv 7$} &impossible &$n\equiv 2,12 \mod 15$&$n\equiv
3,8\mod 15$&impossible&$n\equiv 9,14 \mod 15$\\
 \hline%%%%%%%%%%%%%%%%%%%%%%%%%%%%%%%%%%%%%%%%%%%
 \multicolumn{1}{c|}{$\lt\equiv 8$} &impossible &$n\equiv 4,9 \mod
15$&impossible&$n\equiv 3,8\mod 15$&$n\equiv 2,12 \mod 15$\\
 \hline%%%%%%%%%%%%%%%%%%%%%%%%%%%%%%%%%%%%%%%%%%%
 \multicolumn{1}{c|}{$\lt\equiv 9$} &$n\equiv 0 \mod 5$ &impossible&$n\equiv 2
\mod 5$&$n\equiv 4\mod 5$&$n\equiv 3 \mod 5$\\
 \hline%%%%%%%%%%%%%%%%%%%%%%%%%%%%%%%%%%%%%%%%%%%
 \multicolumn{1}{c|}{$\lt\equiv 10$} &$n\equiv 0 \mod 3$ &$n\equiv 6,11 \mod
15$&$n\equiv 6,11 \mod 15$&$n\equiv 6\mod 15$&$n\equiv 6,11 \mod 15$\\
 \hline%%%%%%%%%%%%%%%%%%%%%%%%%%%%%%%%%%%%%%%%%%%
 \multicolumn{1}{c|}{$\lt\equiv 11$} &impossible  &$n\equiv 3\mod 15$&$n\equiv
9,14 \mod 15$&$n\equiv 2,12\mod 15$&impossible\\
 \hline%%%%%%%%%%%%%%%%%%%%%%%%%%%%%%%%%%%%%%%%%%%
 \multicolumn{1}{c|}{$\lt\equiv 12$} &$n\equiv 0\mod 5$  &$n\equiv 2\mod
5$&$n\equiv 3,13\mod 15$&impossible&$n\equiv 4,9 \mod 15$\\
 \hline%%%%%%%%%%%%%%%%%%%%%%%%%%%%%%%%%%%%%%%%%%%
 \multicolumn{1}{c|}{$\lt\equiv 13$} &impossible  &$n\equiv 9,14\mod
15$&impossible&$n\equiv 3,8\mod 15$&$n\equiv 2,12\mod 15$\\
 \hline%%%%%%%%%%%%%%%%%%%%%%%%%%%%%%%%%%%%%%%%%%%
 \multicolumn{1}{c|}{$\lt\equiv 14$} &impossible  &impossible&$n\equiv 2,12\mod
15$&$n\equiv 9,14\mod 15$&$n\equiv 3 \mod 5$\\
 \hline%%%%%%%%%%%%%%%%%%%%%%%%%%%%%%%%%%%%%%%%%%%
\end{tabular}\label{possiblen}
}

\end{table}
\end{center}

There are at least two other necessary conditions: 

\begin{theorem} Suppose a GDD$(n,2,6;\lambda_1,\lambda_2)$ exists. Then: 
\begin{enumerate}
\item $b\geq\max(2r-\lambda_1,2r-\lambda_2)$
\item $\lambda_2\leq 2\lambda_1(n-1)/n$
\end{enumerate}
\end{theorem}

\begin{proof}
For condition (1), consider the set of blocks containing the points $x$ and $y$.
There are $r$ blocks containing $x$ and $r-\lambda_i$ blocks which contain $y$
but do not contain $x$. So there are at least $2r-\lambda_i$ blocks. For
condition
(2) let $b_6$ be the number of blocks with all 6 points from one group,  $b_5$
be the number of blocks with 5 points from 1 group, and the remaining point from
the other group, $b_4$ be the number of blocks with 4 points from 1 group, and
the remaining 2 points from the other group, and $b_3$ be the number of blocks
with 3 points from each group. Counting the contribution of these blocks towards
the number of pairs of points from the same group in the blocks together gives:
$15b_6+10b_5+7b_4+6b_3=2\lambda_1{n \choose 2}=n(n-1)\lambda_1$. Counting the
pairs of points from different groups gives $5b_5+8b_4+9b_3=n^2\lambda_2$. Thus
we have: 

$-15b_6-5b_5+b_4+3b_3=n^2\lambda_2-n^2\lambda_1+n\lambda_1\leq b_4+3b_3\leq
5b=\\n[\lambda_1(n-1)+\lambda_2n]/3$\\
$\Rightarrow 3n^2\lambda_2-3n^2\lambda_1+3n\lambda_1\leq
n^2\lambda_2+n^2\lambda_1-n\lambda_1$\\
$\Rightarrow 2n^2\lambda_2\leq 4n^2\lambda_1-4n\lambda_1$\\
$\Rightarrow \lambda_2\leq\displaystyle \frac{2(n-1)\lambda_1}{n}$ \end{proof}
 Condition (2) shows that while $\lambda_2\geq \lambda_1$ is possible, we always
have $\lambda_2<2\lambda_1$. We can apply the theorem to assert the following:

 \begin{corollary} The family GDD$(n,2,6;s,2st)$ does not exist for any integers
$s,t>0$. \end{corollary}

In \cite{5gdd14}, Hurd, Mishra and Sarvate proved the following two
results for GDDs with fixed block configuration. We repeat their results
here.

\begin{theorem}[\cite{5gdd14}]\label{evenblocks}
 Suppose a GDD$(n,2,k;\lambda_1,\lambda_2)$ has configuration $(s,t)$. Then the
number of blocks with $s$ points (respectively $t$) from the first group is
equal to the number of blocks with $s$ points (respectively $t$) from the
second group. Consequently, for any $s$ and $t$, the number of blocks $b$ is
necessarily even. 
\end{theorem}

\begin{theorem}[\cite{5gdd14}]\label{betahurd} For any
GDD$(n,2,k;\lambda_1,\lambda_2)$ with
configuration $(s,t)$, the second index is given by
$\lambda_2=\displaystyle\left(\frac{\lambda_1(n-1)}{n}\right)\left(\frac{
k(k-1)-2\beta}{
2\beta}\right)$ where $\beta = \displaystyle{s\choose 2}+{t \choose 2}$.  
\end{theorem}
For the remainder of this paper, we
will refer to the results in this section as the ``necessary'' conditions. 
 %%%%
 %%% gdd's with config 3,3 
 %%%%
%%%%%%%%%%%%%%%%%%%%%%%%%%%%%%%%%%%
 \section{GDDs with Configuration (3,3)}
In this section, we introduce a basic construction for configuration (3,3) GDDs
with specific indices and present the minimal indices for any configuration
(3,3)
GDD($n,2,6;$ $\lo, \lt$). We begin by providing an example of a configuration
(3,3)
GDD where $\lo=4$ and $\lt=5$. \\
 {\bf Example 1}:  GDD(6,2,6;4,5). Let $A=\{0,1,2,3,4,5\}$ and
$B=\{a,b,c,d,e,f\}.$ Then the $b=20$ blocks are: 
\begin{center}
$\{0,1,2,a,b,c\}, \{0,1,2,d,e,f\},  \{0,1,3,a,b,d\},  \{0,1,3,c,e,f\}, 
\{0,2,4,a,c,e\},$\\
 $\{0,2,4,b,d,f\},  \{0,3,5,a,d,f\}, \{0,3,5,b,c,e\}, \{0,4,5,a,e,f\},
\{0,4,5,b,c,d\}, $\\
 $\{1,2,5,b,c,f\}, \{1,2,5,a,e,d\},\{1,3,4,b,d,e\}, \{1,3,4,a,c,e\},
\{1,4,5,b,e,f\}, $\\
 $\{1,4,5,a,c,d\}, \{2,3,4,c,d,e\}, \{2,3,4,a,b,f\}, \{2,3,5,c,d,f\},
\{2,3,5,a,b,e\}$ \end{center}

By applying Theorem \ref{betahurd} to configuration (3,3) GDDs, we get the
following result.

\begin{corollary}\label{l2eq} For any configuration $(3,3)$
GDD$(n,2,6;\lambda_1,\lambda_2)$, we have
$\lambda_2=\displaystyle\frac{3\lo(n-1)}{2n}$.
\end{corollary} 

%%%
%%% A Basic Construction
%%%
\subsection{A Basic Construction for Configuration $(3,3)$}
A balanced incomplete block design BIBD$(v,k,\lambda)$ is a pair $(V,B)$ where
$V$ is a set of points with cardinality $v$ and $B$ is a collection of $b$
$k$-subsets of $V$ called blocks such that each element of $V$ is contained in
exactly $r$ blocks and any 2-subset of $V$ is contained in exactly $\lambda$
blocks. If $k=3$, we may call the design a triple system, and abbreviate
TS$(v,\lambda)$. We use triple systems in the follow construction.  

\begin{theorem}\label{basic} If there exists a TS$(n,\lambda)$ with $b$ blocks
and repetition number $r$, then there
exists a configuration $(3,3)$ GDD$(n,2,6;\lambda b,r^2)$. Further if such a GDD
exists, then there exists a TS$(n,\lambda b)$.  
\end{theorem}

\begin{proof}
Suppose there exists a TS($n,\lambda$). Consider two copies of this triple
system, TS$_1(n,\lambda)$ and TS$_2(n,\lambda)$. Form the complete bipartite
graph $G$ with bipartitions $G_1$ and $G_2$ where $V(G_1)$ is the set of blocks
of TS$_1(n,\lambda)$ and  $V(G_2)$ is the set of blocks of TS$_2(n,\lambda)$.
The blocks of the desired design are the edge set of $G$. Consider a pair of
first associates. They will appear $\lambda$ times in TS$_i(n,\lambda), i=1,2$.
Therefore, in the given construction they will appear together exactly $\lambda
b$ times, where $b$ is the number of blocks in a TS($n,\lambda$). Now consider a
pair of second associates $\{v_1,v_2\}$ where $v_i\in TS_i(n,\lambda)$. Any
point appears exactly $r$ times in a TS($n,\la$), thus the pair $\{v_1,v_2\}$ is
contained in exactly $r^2$ blocks of this design. 

Now suppose a GDD exists with groups $G_1$ and $G_2$.  For each block, remove
the points contained in $G_1$,
and then remove $G_1$.  What remains is a set of blocks of size 3 on $G_2$ which
have the property that any pair
of points occurs in exactly $\lambda b$ blocks.  Thus it is a TS($n, \lambda
b)$.
\end{proof}

The construction given in Theorem \ref{basic} can easily be generalized to any
configuration ($k,k$) GDD. Thus we have the following corollary. 

\begin{corollary}\label{configkk}
If there exists a BIBD$(n,k,\lambda)$ with $b$ blocks and repetition number $r$,
then there exists a configuration $(k,k)$
GDD$(n,2,2k; \lambda b, r^2)$. 
\end{corollary}

%%
%% Minimal Indices
%%%%
\subsection{Minimal Indices}
There exists a TS(7,1), and thus by Theorem \ref{basic} there exists a
GDD(7,2,6;7,9). From Corollary \ref{l2eq},
$\lambda_2=\frac{3\lambda_1(6)}{14}=\frac{9\lambda_1}{7}$, so the construction
given in Theorem \ref{basic} gives a design with the minimum possible indices. 
However,
there also exists a TS(9,1) which means that there exists a GDD(9,2,6;12,16) by
Theorem \ref{basic}. In this case we have that
$\lambda_2=\frac{3\lambda_1(8)}{18}=\frac{4\lambda_1}{3}$. Here the minimum
values for $(\lo,\lt)$ are (3,4). So the construction given in Theorem
\ref{basic} does not give a design with the minimum possible indices. 
In general, Theorem \ref{l2eq} says that for any configuration (3,3) GDD, if for
some value of $n$, the
minimum possible indices are (\lo,\lt), then any other GDD with that
configuration will have the indices ($w$\lo,$w$\lt) for some positive integer
$w$. We can find the minimal indices by using Theorem \ref{l2eq} and by the
equations given in Theorem \ref{thm1}. Any configuration (3,3) GDD with indices
($w$\lo,$w$\lt) can be obtained by taking $w$ copies of the blocks in the
minimal design. Therefore, we focus on constructing configuration (3,3) GDDs
with
indices (\lo, \lt). We may then say that the necessary conditions are sufficient
for the existence of any configuration (3,3) GDD with that $n$.

\begin{table}
\begin{center}
\caption{Summary of Minimal Indices for Configuration $(3,3)$}
\begin{tabular}{l|r|r}
 \hline%%0
 \multicolumn{1}{c|}{$n$} &$\lambda_1$ &$\lambda_2$  \\%%1
\hline
\multicolumn{1}{c|}{$n\equiv 0\mod 6$} & $2n/3$ &$(n-1)$\\%%1
 \hline
\multicolumn{1}{c|}{$n\equiv 1\mod 6 $} &  $n$ & $3(n-1)/2$\\%%2
 \hline
\multicolumn{1}{c|}{$n\equiv 2 \mod 6$} &$6n$&$9(n-1)$ \\%%3
 \hline
\multicolumn{1}{c|}{$n\equiv 3 \mod 6$} &  $n/3$
&$(n-1)/2$\\%%4
\hline
\multicolumn{1}{c|}{$n\equiv 4 \mod 6$} &$2n$ & $3(n-1)$ \\%%5
\hline
\multicolumn{1}{c|}{$n\equiv 5 \mod 6$} & $3n$ &$9(n-1)/2$\\
\hline
\end{tabular}\label{33table}
\end{center}
\end{table}
\begin{theorem}\label{lambdaTable33}
The minimal indices $(\lo,\lt)$ for any configuration $(3,3)$
GDD$($$n$,$2$, $6$;$\lo$,$\lt$$)$ are summarized in Table \ref{33table}. 

\end{theorem}

\begin{proof}

We know that $\lt=\frac{3\lo (n-1)}{2n}$ from Theorem \ref{l2eq}. If $n\equiv
0\mod 3$ and
$n\equiv 1\mod 2$, then $n\equiv 3 \mod 6$. Thus $\lo$ is a multiple of $n/3$
and $\lt$ is a multiple of $(n-1)/2$. If $n\equiv 0\mod 3$ and $n\equiv 0\mod
2$, then $n\equiv 0 \mod 6$, so $\lo$ is a multiple of $2n/3$ and $\lt$ is a
multiple of $(n-1)$. If $n\equiv 1\mod 3$ and $n\equiv 1\mod 2$, $n\equiv 1 \mod
6$, implying $\lo$ is a multiple of $n$ and $\lt$ is a multiple of $3(n-1)/2$.
If $n\equiv 1\mod 3$ and $n\equiv 0\mod 2$, $n\equiv 4 \mod 6$, and  $\lo$ is a
multiple of $2n$ and $\lt$ is a multiple of $3(n-1)$. If $n\equiv 2\mod 3$ and
$n\equiv 1\mod 2$, then $n\equiv 5 \mod 6$. This implies that $\lo$ is a
multiple of $n$ and $\lt$ is a multiple of $3(n-1)/2$. However, if we take these
values to be the minimal indices, these number of blocks given by Theorem
\ref{thm1} would not be integer valued. The smallest values for $(\lo,\lt)$ that
give integer values for $b$ are $(\lo, \lt)=(3n,\frac{9}{2}(n-1))$. Finally
consider when $n\equiv 2 \mod 3$ and $n\equiv 0 \mod 2$. Then $n\equiv 2\mod 6$,
which means that $\lo$ is a multiple of $2n$ and $\lt$ is a multiple of
$3(n-1)$. If we take these values to be the minimal indices, these number of
blocks given by Theorem \ref{thm1} would not be integer valued so the smallest
values for $(\lo,\lt)$ that give integer values for $b$ are $(\lo,
\lt)=(6n,9(n-1))$.

\end{proof}

\section{Constructing Configuration (3,3) GDDs}
In this section, we give a similar construction to the one given in Theorem
\ref{basic} based
on $\alpha$-resolvable triple systems. We then show that this construction
produces
designs with minimal indices for all configuration (3,3) GDDs with block size 6
and 2 groups.  

A set of blocks in a design is called a parallel class if it partitions the
point set. A partition of the blocks of a design into parallel classes is a
resolution, and such a design is called resolvable. An $\alpha$-parallel class
in a design is a set of blocks which contain every point of the design exactly
$\alpha$ times. A design that can be resolved into $\alpha$-parallel classes is
called $\alpha$-resolvable. We may abbreviate an $\alpha$-resolvable design as
an
$\alpha$-RBIBD$(n,k,\lambda)$. If $\alpha=1$ then we abbreviate
RBIBD$(n,k,\lambda)$. 

The necessary conditions for the existence of a $\alpha$-RBIBD$(n,k\lambda)$
were given by Jungnickle, Mullin and Vanstone in \cite{alphres3}.

\begin{theorem}[\cite{alphres3}]\label{alphares}
The necessary conditions for the existence of an $\alpha$-resolvable
BIBD$(n,k,\lambda)$ are 
\begin{enumerate}
\item $\lambda(n-1)\equiv 0 \mod{(k-1)\alpha}$
\item $\lambda n(n-1)\equiv 0 \mod{k(k-1)}$
\item  $\alpha n \equiv 0 \mod k$ 
\end{enumerate}
\end{theorem}
In the same paper, they also showed that these conditions were sufficient when
$k=3$. 
\begin{lemma}[\cite{alphres3}]\label{alphares3}
The necessary conditions for the existence of an $\alpha$-resolvable
BIBD$(n,$ $3,$ $\lambda)$ are sufficient, except for $n=6, \alpha=1$ and
$\lambda\equiv 2 \mod 4$. 
\end{lemma}
Vasiga, Furino and Ling \cite{alphares4} showed that the necessary conditions
are
sufficient for $k=4$. 
\begin{lemma}[\cite{alphares4}]\label{alphares4}
The necessary conditions for the existence of an $\alpha$-resolvable
BIBD$(n,$ $4,$ $\lambda)$ are sufficient, with the exception of
$(\alpha,n,\lambda)=(2,10,2)$. 
\end{lemma}

We use $\alpha$-resolvable designs to obtain the following result. 
%%%%%
%%%%%
%%% master lemma
%%%%%
%%%%%
\begin{lemma}\label{masterlemma}
Suppose there exists an $\alpha$-resolvable TS$(n,\la)$ with $s$
$\alpha$-parallel classes, where each parallel class contains $t$ blocks. Then
there exists a configuration $(3,3)$ GDD$(n,2,6;\la t, \alpha^2 s)$. 
\end{lemma}

\begin{proof}
For $i=1,2$, let $D_i$ be an $\alpha$-resolvable TS$(n,\la)$. Resolve the blocks
of $D_i$ into $\alpha$-parallel classes $C_1^i, C_2^i, \dots C^i_s$. Construct a
graph $G$ in the following manner.  For $j=1, 2, \dots, s$, create the complete
bipartite graph $G_j$ with bipartitions $G_j^1$ and $G_j^2$ where $V(G_j^1)$ are
the blocks of $C_j^1$ and $V(G_j^2)$ are the blocks of $C_j^2$. Let
$G=\bigcup^s_{j=1}G_j$. The edge set of $G$ will form the blocks of the desired
design. 

Consider a pair of first associates. It will appear in exactly \la\; blocks of
$D_i$. Therefore, in the given construction, it will appear in $\la t$ blocks of
size 6. Now consider a pair of second associates $\{v_1, v_2\}$ where $v_1\in
D_1$ and $v_2\in D_2$. Here $v_1$ will be matched with $v_2$ exactly $\alpha^2$
times per $\alpha$-parallel class, thus $\lambda_2=\alpha^2 s$.\end{proof}

We now consider values of $n\mod 6$ and apply Lemma \ref{masterlemma} in each
case to obtain the desired configuration (3,3) GDD with minimal indices (\lo,
\lt). 
%%%%%
%%%%%
%%% n congruent to 3
%%%%%
%%%%%
\begin{theorem} The necessary conditions are sufficient for the existence of
a configuration $(3,3)$ GDD
$(n,2,6; \frac{n}{3},\frac{n-1}{2})$ when $n\equiv 3\mod{6}$. 
\end{theorem}

\begin{proof}
Let $n\equiv 3 \mod 6$. Then by Lemma \ref{alphares3} there exists a
1-resolvable TS($n,1$) with $\frac{n-1}{2}$ parallel classes, each containing
$\frac{n}{3}$ blocks. By applying the construction in Lemma \ref{masterlemma} we
obtain a GDD with indices $(\lo,\lt)=(\frac{n}{3}, \frac{n-1}{2})$, which are
 the minimal indices given in Theorem \ref{lambdaTable33}.\end{proof}
%%%%%
%%%%%
%%% n congruent to 1
%%%%%
%%%%%
\begin{theorem} The necessary conditions are sufficient for the existence of GDD
$(n,2,6;n,\frac{3}{2}(n-1))$ when $n\equiv 1\mod{6}$ with configuration
$(3,3)$. 
\end{theorem}

\begin{proof}
Let $n\equiv 1 \mod 6$. By Lemma \ref{alphares3} there exists a 3-resolvable
TS($n,1$) with $\frac{n-1}{6}$ 3-parallel classes, each containing $n$ blocks.
If
we apply the construction in Lemma \ref{masterlemma}, we obtain a GDD with
minimal indices $(\lo,\lt)=(n, \frac{3(n-1)}{2})$. \end{proof}

\begin{theorem} The necessary conditions are sufficient for the existence of GDD
$(n,2,6;6n,9(n-1))$ when $n\equiv 2\mod{6}$ with configuration $(3,3)$. 
\end{theorem}
%%%%%
%%%%%
%%% n congruent to 2
%%%%%
%%%%%

\begin{proof}
Let $n\equiv 2 \mod 6$. Then by Lemma \ref{alphares3} there exists a
3-resolvable TS($n,6$) with $(n-1)$ 3-parallel classes, each containing $n$
blocks. Applying Lemma \ref{masterlemma} yields a GDD with minimal indices
$(\lo,\lt)=(6n, 9(n-1))$. \end{proof}

\begin{theorem} The necessary conditions are sufficient for the existence of GDD
$(n,2,6;2n,3(n-1))$ when $n\equiv 4\mod{6}$ with configuration $(3,3)$. 
\end{theorem}
%%%%%
%%%%%
%%% n congruent to 4
%%%%%
%%%%%
\begin{proof}
Let $n\equiv 4 \mod 6$. By Lemma \ref{alphares3}, there exists a 3-resolvable
TS($n,2$) with $\frac{n-1}{3}$ 3-parallel classes each containing $n$ blocks. We
may apply Lemma \ref{masterlemma} to obtain a GDD with minimal indices
$(\lo,\lt)=(2n, 3(n-1))$. \end{proof}
%%%%%
%%%%%
%%% n congruent to 5
%%%%%
%%%%%
\begin{theorem} The necessary conditions are sufficient for the existence of GDD
$(n,2,6;3n,\frac{9}{2}(n-1))$ when $n\equiv 5\mod{6}$ with configuration
$(3,3)$. 
\end{theorem}

\begin{proof}
Let $n\equiv 5 \mod 6$. Then by Lemma \ref{alphares3} there exists a
3-resolvable TS($n,3$) with $\frac{n-1}{2}$ 3-parallel classes, each containing
$n$ blocks. We may apply Lemma \ref{masterlemma} to obtain a GDD with minimal
indices $(\lo,\lt)=(3n, \frac{9(n-1)}{2})$. \end{proof}

\begin{theorem}
The necessary conditions are sufficient for the existence of
GDD$(n$, $2$, $6;$ $\frac{2}{3}n$, $n-1)$ for $n\equiv 0 \mod 6$ with
configuration
$(3,3)$. 
\end{theorem}
\begin{proof}
Let $n\equiv 0\mod 6$ with $n\geq 12$. Then by Lemma \ref{alphares3} there
exists a 1-resolvable TS($n,2$) with $n-1$ parallel classes, each containing
$\frac{n}{3}$ blocks. If we apply the construction given in Lemma
\ref{masterlemma} we obtain a GDD with minimal indices $(\lo,\lt)=(\frac{2n}{3},
n-1)$. If $n=6$, we may not use the construction described in Lemma
\ref{alphares3}. However if $n=6$, the minimal indices $(\lo, \lt)=(4,5)$ and
Example 1 gives a GDD(6,2,6;4,5). \end{proof}

Since we have given a construction for all possible values of $n\mod 6$, we may
give the following result.

\begin{theorem}
The necessary conditions are sufficient for the existence of all
configuration $(3,3)$ GDD$(n,2,6;\lo,\lt)$ with minimal indices. 
\end{theorem}  
%%%%%%%%%%%%%%%%%%%%%%%%%%%%%%%%%%%
%%%%%%%%%%%%%%%%%%%%%%%%%%%%%5%%%%%
%CONFIG 4 2
%%%%%%%%%%%%%%%%%%%%%%%%%%%%%%%%%%
\section{GDDs with Configuration (2,4)}
In this section we present the minimal indices for any configuration (2,4)
GDD$(n,$ $2,6;\lo,\lt)$. 
By Theorem \ref{betahurd} we have the following relation between $\lambda_1$
and $\lambda_2$ for any configuration $(2,4)$ GDD. 

\begin{theorem}\label{fourtwolambdas}
 For any configuration $(2,4)$ GDD$(n,2,6;\lambda_1,\lambda_2)$ we have
$\lambda_2=\frac{8\lambda_1(n-1)}{7n}$. 
\end{theorem}
For any configuration (2,4) GDD if for some value of $n$, the minimum
possible indices are $(\lambda_1,\lambda_2)$, then any other GDD with that
configuration will have the indices $(w\lambda_1,w\lambda_2)$ for some positive
integer $w$. We may find the minimum indices by using the equation in Theorem
\ref{fourtwolambdas}, the equations in Theorem \ref{thm1}, and the condition in
Theorem
\ref{evenblocks}. As in the case with configuration $(3,3)$, we focus on
constructing GDDs with minimal indices since we may then say the necessary
conditions are sufficient for the existence of any configuration (2,4) GDD with
that $n$.

\begin{center}
\begin{table}

\caption{Summary of Minimal Indices for Configuration $(2,4)$}
\begin{tabular}{l|r|r}

 \hline%%0
 \multicolumn{1}{c|}{$n$} &$\lambda_1$ &$\lambda_2$  \\%%1
%\multicolumn{1}{c|}{} &  & \\
\hline
\multicolumn{1}{c|}{$n\equiv 0,16,24,32,40,48\mod 56$ } &
$7n/8$ &$n-1$\\%%1
%\multicolumn{1}{c|}{} &  & \\
 \hline
\multicolumn{1}{c|}{$n\equiv 2,6,10,14,18,26,30,34,38,42,46,54 \mod 56 $} & 
$7n/2$ &
$4(n-1)$\\%%2
%\multicolumn{1}{c|}{} &  & \\
 \hline
\multicolumn{1}{c|}{$n\equiv 4,12,20,28,44,52 \mod 56$} &$7n/4$&$2(n-1)$
\\%%3
%\multicolumn{1}{c|}{} &  & \\
 \hline
\multicolumn{1}{c|}{$n\equiv 8 \mod 56$} & $n/8$ &$(n-1)/7$\\%%4
%\multicolumn{1}{c|}{} &  & \\
\hline
\multicolumn{1}{c|}{$n\equiv 22,50 \mod 56$} & $n/2$
&$4(n-1)/7$\\%%4
%\multicolumn{1}{c|}{} &  & \\
\hline
\multicolumn{1}{c|}{$n\equiv 36 \mod 56$} & $n/4$
&$2(n-1)/7$\\%%4
%\multicolumn{1}{c|}{} &  & \\
\hline
\multicolumn{1}{c|}{$n\equiv
3,5,7,9,11,13,17,19,21,23,25,27,31,$} &  & \\%%4
\multicolumn{1}{c|}{$33,35,37,39,41,45,47,49,51,53,55 \mod
56$} &$7n$ & $8(n-1)$
\\
%\multicolumn{1}{c|}{} &  & \\
\hline
\multicolumn{1}{c|}{$n\equiv 1,15,29,43 \mod 56$} &$n$ & $8(n-1)/7$
\\%%5
%\multicolumn{1}{c|}{} &  & \\
\hline
\end{tabular}\label{421table}
\end{table}
\end{center}

\begin{theorem}\label{lambdaTable42}
The minimal indices $(\lo,\lt)$ for any configuration $(2,4)$
GDD$(n,2, $$6;$ $\lambda_1,$ $\lambda_2)$ are summarized in Table
\ref{421table}. 

\end{theorem}
\begin{proof}
By Theorem \ref{fourtwolambdas}, we know that $\lt=\frac{8\lo(n-1)}{7n}$. If
$n\not\equiv 1\mod 7$ and $n$ is odd, then this implies that $n\equiv
3,5,7,9,11,13\mod 14$. Thus $\lo$ is a multiple of $7n$ and $\lt$ is a multiple
of $8(n-1)$. If $n\equiv 1\mod 7$ and $n$ is odd, then $n\equiv 1 \mod 14$. In
this case, $\lo$ must be a multiple of $n$ and $\lt$ a multiple of $(8/7)(n-1)$.
If $n\not \equiv 1\mod 7$ and $n\equiv 0\mod 8$, we have that $n\equiv
0,16,24,32,40,48\mod 56$, so $\lo$ is a multiple of $7n/8$ and $\lt$ is a
multiple of $n-1$.  If $n\not \equiv 1\mod 7$ and $n\equiv 2\mod 8$ then
$n\equiv 2,10,18,26,34,42\mod 56$ implying $\lo$ is a multiple of $7n/2$ and
$\lt$ is a multiple of $4(n-1)$.  If $n\not \equiv 1\mod 7$ and $n\equiv 4\mod
8$, $n\equiv 4,12,20,28,44,52\mod 56$. Then $\lo$ is a multiple of $7n/4$ and
$\lt$ is a multiple of $2(n-1)$.  If $n\not \equiv 1\mod 7$ and $n\equiv 6\mod
8$, $n\equiv 6,14,30,38,46,54\mod 56$, then $\lo$ is a multiple of $7n/2$ and
$\lt$ is a multiple of $4(n-1)$.  If $n\equiv 1\mod 7$ and $n\equiv 0\mod 8$, we
have that $n\equiv 8\mod 56$. Here, it follows that $\lo$ is a multiple of $n/8$
and $\lt$ is a multiple of $(n-1)/7$. If $n\equiv 1\mod 7$ and $n\equiv 2\mod
8$, we have that $n\equiv 50\mod 56$. Here, it follows that $\lo$ is a multiple
of $n/2$ and $\lt$ is a multiple of $4(n-1)/7$. If $n\equiv 1\mod 7$ and
$n\equiv 4\mod 8$, we have that $n\equiv 36\mod 56$. Here, it follows that $\lo$
is a multiple of $n/4$ and $\lt$ is a multiple of $2(n-1)/7$. If $n\equiv 1\mod
7$ and $n\equiv 6\mod 8$, $n\equiv 22\mod 56$, and it follows $\lo$ is a
multiple of $n/2$ and $\lt$ is a multiple of $4(n-1)/7$.
\end{proof}

%%%%%%%%%%%%%%%%%%%%%%%%%%%%%%%%%%%
%%%%%%%%%%%%%%%%%%%%%%%%%%%%%5%%%%%
%CONFIG 4 2 Constructions
%%%%%%%%%%%%%%%%%%%%%%%%%%%%%%%%%%

\section{Constructing $(2,4)$ GDD$(n,2,6;\lo,\lt)$}
We use the Theorem \ref{lambdaTable42} and Lemma \ref{alphares4} to construct
configuration (2,4) GDDs with minimal indices, when possible. We begin with a
general
construction.

\begin{lemma}\label{masterlemma24} 
If there exists an $\alpha$-resolvable BIBD$(n,4,\lambda)$ with $n$ even and
$\lambda=3\alpha$, then there exists a
configuration $(2,4)$ GDD$(n,2,6;\frac{n}{2}(\lambda+\frac{\alpha}{2}),
2\alpha(n-1))$. 
\end{lemma}
\begin{proof}
 Let the two groups be $A=\{1,2,\dots,n\}$, and $A'=\{1',2',\dots,n'\}$. Let $D$
be an $\alpha$-resolvable BIBD$(n,4,\lambda)$ on the point set of $A$. Let $F$
be a 1-factorization of $K_n$ on the point set of $A'$. Resolve
the blocks into $\alpha$ parallel classes. There will be
$\lambda(n-1)/3\alpha=(n-1)$ 
classes with $(n\alpha)/4$ blocks in each class. Construct a graph $G$ in the
following manner. For $j=1,2,\dots, (n-1)$, create the complete
bipartite graph $G_j$ with bipartitions $G_j^1$ and $G_j^2$ where $V(G_j^1)$
are the blocks of an $\alpha$ parallel class and $V(G_j^2)$ are a 1-factor of
$K_n$. If we switch $A$ with $A'$ and repeat the construction, we obtain all
desired blocks. 

Consider a pair of first associates, $\{x,y\}\in A$. It will appear exactly
$\lambda$ times in
$D$. Therefore in the given construction, it will appear $n\lambda/2$
times when in the first part of the construction. This pair will appear an
additional
$n\alpha/4$ times when the second part of the construction. Thus
$\lo=\frac{n}{2}(\lambda+\frac{\alpha}{2})$. Now consider a pair of second
associates
$\{x,y'\}$, where $x\in A$ and $y'\in A'$. Here $x$ will appear with $y'$
exactly $\alpha(n-1)$ times in both
parts of the construction, so $\lt=2\alpha(n-1)$. 
\end{proof}

We use the above construction to obtain the following results:

\begin{corollary}
 Let $n\equiv 2,6,10,14,18,26,30,34,38,42,46,54 \mod 56$. Then the necessary
conditions are sufficient 
for the existence of a configuration $(2,4)$
GDD$(n,$ $2,6;$ $\frac{7n}{2},4 (n-1))$. 
\end{corollary}
\begin{proof}
Let $n$ be assumed as above. By Lemma \ref{alphares4}, there exists a
2-resolvable BIBD$(n,4,6)$. Apply Lemma \ref{masterlemma24}
to obtain a GDD with minimal indices $(\lambda_1,\lambda_2)=(\frac{7n}{2},4
(n-1))$. 
\end{proof}

\begin{corollary}

 Let $n\equiv 4,12,20,28,44,52 \mod 56$. Then the necessary
conditions are sufficient 
for the existence of a configuration $(2,4)$
GDD$(n,2,6$; $\frac{7n}{4},2 (n-1))$. 
\end{corollary}
\begin{proof}
Let $n$ be assumed as above. By Lemma \ref{alphares4}, there exists a
resolvable BIBD$(n,4,3)$. So we may apply Lemma \ref{masterlemma24}
to obtain a GDD with minimal indices $(\lambda_1,\lambda_2)=(\frac{7n}{4},2
(n-1))$. 
\end{proof}
 We define a near-minimal GDD as a GDD which has indices exactly twice the
minimal size. 

\begin{corollary}
 If $n\equiv 0,8\mod 24$ then there exists a near minimal configuration $(2,4)$
GDD$(n,2,6;\frac{7n}{4},2 (n-1))$. 
\end{corollary}
\begin{proof}
Let $n$ be assumed as above. By Lemma \ref{alphares4}, there exists a
resolvable BIBD$(n,4,3)$. Apply Lemma \ref{masterlemma24}
to obtain a near-minimal GDD with indices
$(\frac{7n}{4},2(n-1))$. 
\end{proof}

The above construction gives near-minimal GDDs for $n\equiv 0, 8\mod 24$. The
next theorem shows that for $n=8$, the minimal indices can not be obtained. 

\begin{theorem}
 There does not a exist a configuration $(2,4)$ GDD$(8,2,6;1,1)$. 
\end{theorem}
\begin{proof}
 Assume such a design exists with groups $A$ and $B$. Then it would have 8
blocks and every point would
appear 3 times. Consider a point in the design, $x$ and let its first
associates be $\{1,2,3,4,5,6,7\}$. Suppose $x$ appears
in 3 blocks which intersect $A$ in 4 points, and $x\in A$ in each of these
blocks. Then because there are only 7 other points,
there must be a repeated pair in one of these blocks. However, we
assumed $\lambda_1=1$, so this is not possible. Now suppose $x$ appears in 2
blocks which intersect $A$ in 4 points and $x\in A$ in those blocks.  Then $x$
must also appear in a block which intersects $B$ in 2 points and $x\in B$.
Let the two partial blocks containing $x\in A$ be $\{x,1,2,3\}$ and
$\{x,4,5,6\}$.
Without loss, assume
the last partial block containing $x$ also contains $1$, and $1\in A$. The part
of this block which intersects $A$ may not contain $x,2,3$, and we cannot
repeat pairs, so $1$ must be in a partial block with $\{4,7\}$. However, there
is
no additional first associate available to complete this block. Finally assume
$x$ appears in one block which
intersects $A$ in 4 points and $x\in A$. Without loss, we may assume the partial
block containing $x\in A$ be $\{x,1,2,3\}$. Then $x$ appears in 2 blocks which
intersect $B$ in 2 points, and $x\in B$. One of these blocks must contain the
pair $\{x,1\}$ where $1\in A$ and the other block must contain the pair
$\{x,2\}$ where $2\in A$. However, we have no way to cover the pair $\{x,3\}$
where $x\in A$ and $3\in B$ or $x\in B$ or $3\in A$. Thus this design cannot
exist.
\end{proof}

We use a slightly different construction for $n\equiv 16 \mod 24$. 

\begin{theorem}
If $n\equiv 16 \mod 24$ then the necessary conditions are sufficient for the
existence of a configuration $(2,4)$ GDD$(n,2,6;7n/8,(n-1))$.  

\end{theorem}

\begin{proof}
 Let $n\equiv 16\mod 24$, and let $A=\{1,2,\dots,n\}$ and
$A'=\{1',2',\dots,n'\}$ be
the point set for the two groups in the desired design. By Theorem
\ref{alphares4}, there exists a RBIBD$(n,4,1)$. Let $D$ be such a design with
point set $A$. Resolve the blocks of $D$ into parallel classes, $C_1,\dots,
C_{(n-1)/3}$. There will be  $n/4$ blocks in each
parallel class. We construct a 1-factorization of $K_n$ on the point set of
$A'$. On each parallel class $C_j,j=1,2,\dots,(n-1)/3$, decompose the blocks of
$C_j$ into three 1-factors as follows. For each block $\{a,b,c,d\}\in C_j$ we
let
$\{\{a',b'\}, \{c',d'\}\}\in F_{j,1}$, $\{\{a',c'\}, \{b',d'\}\}\in F_{j,2}$ and
$\{\{a',d'\},
\{b',c'\}\}\in F_{j,3}$. 

Now construct a graph $G$ in the following manner. For $j=1,2,\dots,(n-1)/3$,
construct the complete bipartite graph $G_{j,1}$ with bipartitions $G^1_{j,1}$
and
$G^2_{j,1}$ where $V(G^1_{j,1})$ are, without loss, the first $n/8$ blocks of
parallel class
$C_j$ and $V(G^2_{j,1})$ are the 1-factor $F_{j,1}$. Also, create the complete
bipartite graph $G_{j,2}$ with bipartitions $G^1_{j,2}$ and
$G^2_{j,2}$ where $V(G^1_{j,2})$ are the last, without loss, $n/8$ blocks of
parallel class
$C_j$ and $V(G^2_{j,2})$ are the 1-factor $F_{j,2}$. Construct the complete
bipartite graph $G_{j,3}$ with bipartitions $G^1_{j,3}$ and
$G^2_{j,3}$ where $V(G^1_{j,3})$ are the first $n/8$ blocks of parallel class
$C_j$ and $V(G^2_{j,3})$ are the edges of $F_{j,3}$ which were obtained from
the first $n/8$ blocks of $C_j$. Finally construct the complete bipartite graph 
$G_{j,4}$ with bipartitions
$G^1_{j,4}$ and
$G^2_{j,4}$ where $V(G^1_{j,4})$ are the last $n/8$ blocks of parallel class
$C_j$ and $V(G^2_{j,4})$ are the edges of $F_{j,3}$ which are obtained from
the last $n/8$ blocks of $C_j$. If we take the union of all these bipartite
graphs, then we obtain half of the blocks of size 6 in the GDD. 

To obtain the other half, we switch the roles of $A$ and $A'$ in the design and
the 1-factorization. We construct a graph $H$ on the vertex set $A,A'$ in a
similar manner to $G$. For $j=1,2,\dots,(n-1)/3$,
construct the complete bipartite graph $H_{j,1}$ with bipartitions $H^1_{j,1}$
and
$H^2_{j,1}$ where $V(H^1_{j,1})$ are the last $n/8$ blocks of parallel class
$C_j$ and $V(H^2_{j,1})$ are the 1-factor $F_{j,1}$. Also, construct the
complete
bipartite graph $H_{j,2}$ with bipartitions $H^1_{j,2}$ and
$H^2_{j,2}$ where $V(H^1_{j,2})$ are the first $n/8$ blocks of parallel class
$C_j$ and $V(H^2_{j,2})$ are the 1-factor $F_{j,2}$. Construct the complete
bipartite graph $H_{j,3}$ with bipartitions $H^1_{j,3}$ and
$H^2_{j,3}$ where $V(H^1_{j,3})$ are the first $n/8$ blocks of parallel class
$C_j$ and $V(H^2_{j,3})$ are the edges of $F_{j,3}$ which were obtained from
the last $n/8$ blocks of $C_j$. Finally construct the complete bipartite graph 
$H_{j,4}$ with bipartitions
$H^1_{j,4}$ and
$H^2_{j,4}$ where $V(H^1_{j,4})$ are the last $n/8$ blocks of parallel class
$C_j$ and $V(H^2_{j,4})$ are the edges of $F_{j,3}$ which are obtained from
the first $n/8$ blocks of $C_j$. If we take the union of all these bipartite
graphs, then we obtain the other half of the blocks of size 6 in the GDD. 

Consider a pair of first associates. In the first part of the construction, when
$\{x,y\}\in A$ appears in the BIBD, it will
appear exactly once. Thus in the construction, it will be in a block of size 6
exactly $n/2+n/4=3n/4$ times. In the second part of the construction when
$\{x,y\}$ is in the role of a 1-factor, it will appear $n/8$ times. Thus
$\lo=7n/8$. Now consider a pair of second associates,
$\{x,y'\}$ where $x\in A$ and $y'\in A'$. Without loss, we may assume
$\{x,y'\}\in C_j$ for some $j$. In part one of the construction, there are 4
cases to consider. Each point is either in the first $n/8$ blocks of $C_j$ or
in the last $n/8$ blocks of $C_j$. Let $C_{j,1}$ denote the first $n/8$ blocks
of
$C_j$ and $C_{j,2}$ denote the last $n/8$ blocks of $C_j$. Suppose $x\in
C_{j,1}$ and $y'\in C_{j,1}$. Then in the construction, $\{x,y'\}$ appears
twice. If $x\in C_{j,1}$ and $y'\in C_{j,2} $, then $\{x,y'\}$ appears once. If
$x\in C_{j,2}$ and $y'\in C_{j,1} $, then $\{x,y'\}$ appears once and if $x\in
C_{j,2}$ and $y'\in C_{j,2} $, then $\{x,y'\}$ appears twice. In the second
part of the construction when we reverse the roles, if $x\in C_{j,1}$ and $y'\in
C_{j,1}$, then $\{x,y'\}$ appears once.  If $x\in C_{j,1}$ and $y'\in C_{j,2}$,
then
$\{x,y'\}$ appears twice.  If $x\in C_{j,2}$ and $y'\in C_{j,1}$, then
$\{x,y'\}$ appears twice, and if $x\in C_{j,2}$ and $y'\in C_{j,2}$, then
$\{x,y'\}$ appears once.  Thus for each parallel class, each pair $\{x,y'\}$
appears a total of 3 times. Thus each pair of second associates will appear a
total of $3(n-1/3)=n-1$ times in the construction.
\end{proof}

\begin{theorem}
 Let $n\equiv 3,5,7,9,11,13 \mod{14}$. Then the necessary conditions are
sufficient 
for the existence of a configuration $(2,4)$ GDD$(n,2,6;7n,8(n-1))$. 
\end{theorem}
\begin{proof}
  Let the two groups be $A=\{1,2,\dots,n\}$ and $A'=\{1',2',\dots,n'\}$. By
Lemma \ref{alphares4}, there exists a 4-resolvable BIBD($n,4,6$). Let $D$ be
such a
design
with point set $A$. Resolve the blocks of $D$ into 4-parallel classes. There
will be $(n-1)/2$ classes with $n$ blocks in each class. Construct a graph $G$
in the following manner. For $j=1,2,\dots,(n-1)/2$ create the complete bipartite
graph $G_j$ with bipartitions $G_j^1$ and $G_j^2$ where $V(G_j^1)$ are the
blocks of a 4-parallel class and $V(G_j^2)$ are the pairs obtained by
developing $\{0',j'\}\mod n$. If we switch $A$ with $A'$ and repeat the same
construction, we obtain
all desired blocks. 
 
 Consider a pair of first associates, $\{x,y\}\in A$. It will appear exactly 6
times in $D$.
Therefore, in the given construction, it will appear $6n$ times in the first
part of the construction. This pair will appear an additional $n$ times when 
in the second part. Thus $\lambda_1=7n$. Now consider a pair of second
associates $\{x,y'\}$
where $x\in A$ and $y'\in A'$. Here $x$ will
be matched with $y'$ exactly $4(n-1)$ times, in each part of the construction,
and thus $\lambda_2=8(n-1)$. 

\end{proof}

If $n\equiv 1,15,29,43\mod 56$, then the above construction gives a GDD with 7
times the minimal indices. However, the following construction gives a
configuration $(2,4)$ GDD$(15,2,6;15,16)$ with minimum possible indices. 

%%%%%%%%5

\begin{theorem}
The necessary conditions are sufficient for the existence of a
configuration $(2,4)$ GDD$(15,2,6;15,16)$. 
\end{theorem}

\begin{proof}
 By Lemma \ref{alphares4}, there exists a RBIBD($16,4,1$).  It has 5 parallel
classes with 4 blocks in
each class.  Let $X=\{\infty, 0,1,2, ..., 14\}$ be the
points in the RBIBD($16,4,1$).  Because $\infty$ appears with every other point
exactly once, the blocks of the form $\{\infty, x,y,z\}$ form a partition the
set $X\backslash\{\infty\}$.  Each block is in one of the 5 parallel classes. 
For each block $\{\infty, x,y,z\}$, form the pairs $\{x,y\}, \{x,z\}, \{y,z\}$.
Let the two groups be $A=\{0,1,...,14\}$ and $A'=\{0',1', ..., 14'\}$. 
For $j=1,2,3,4,5$, create the complete bipartite graph G$_j$ with bipartitions
$G_{j_1}$ and $G_{j_2}$ where $V(G_{j_1})$ are the blocks of parallel class $j$
except the block containing $\infty$, and $V(G_{j_2})$ are the 15 pairs obtained
from the blocks containing $\infty$. This gives us half of the desired blocks.
To get the rest of the blocks repeat the construction with $V(G_{j_1})$
as the 15 pairs and $V(G_{j_2})$ as the blocks of $PC_j$. 

Consider a pair of first associates, $\{x,y\}\in A$.  If $\{x,y\}$ was in a
block
with $\infty$ in the RBIBD, then it appears exactly 0 times in the first part of
the construction
and 15
times in the second part. If $\{x,y\}$ was not in a block with $\infty$ in the
RBIBD,
then it appears exactly 15 times in the first part and 0 times in the second
part. Therefore,
each pair of first associates appears $\lambda_1=15$ times.  Now consider a pair
of second associates $\{x,y'\}$ where $x\in A$ and $y'\in A'$.  In the first
part, $x$ is in 4 of the blocks
and $y'$
is in 2 of the blocks, so $\{x,y'\}$ is in 8 blocks.  In the second part, $x$ is
in 2
blocks and $y'$ is in 4 blocks, so $\{x,y'\}$ is again in 8 blocks.  Thus,
$\lambda_2=16$.
\end{proof}
\subsection{Summary of Minimality}

\begin{center}
\begin{table}[here]
\caption{Summary of Constructions and Minimality for Configuration $(2,4)$}
\scalebox{0.8}{
\begin{tabular}{l|r|r|| c}

 \hline%%0
 \multicolumn{1}{c|}{$n$} &$\lambda_1$ &$\lambda_2$  \\%%1
\hline
\multicolumn{1}{c|}{$n\equiv 0,16,24,32,40,48\mod 56$ and $n\equiv 16 \mod 24$ }
& $7n/8$ &$n-1$ & minimal\\%%1
 \hline
\multicolumn{1}{c|}{$n\equiv 0,16,24,32,40,48\mod 56$ and $n\equiv 0,8 \mod 24$
}
& $7n/8$ &$n-1$ & near minimal\\%%1
 \hline
\multicolumn{1}{c|}{$n\equiv 2,10,18,26,34,42,6,14,30,38,46,54 \mod 56 $} & 
$7n/2$ &
$4(n-1)$ & minimal\\%%2
 \hline
\multicolumn{1}{c|}{$n\equiv 4,12,20,28,44,52 \mod 56$} &$7n/4$&$2(n-1)$ &
minimal \\%%3
 \hline
\multicolumn{1}{c|}{$n\equiv 8 \mod 56$ and $n\equiv 16 \mod 24$} &
$n/8$ &$(n-1)/7$& 7 times the minimal\\%%4
\hline
\multicolumn{1}{c|}{$n\equiv 8 \mod 56$ and $n\equiv 0,8 \mod 24$} &
$n/8$ &$(n-1)/7$& 14 times the minimal\\%%4
\hline
\multicolumn{1}{c|}{$n\equiv 22,50 \mod 56$} & $n/2$
&$4(n-1)/7$& 7 times the minimal\\%%4
\hline
\multicolumn{1}{c|}{$n\equiv 36 \mod 56$} & $n/4$ &$2(n-1)/7$&
 7 times the minimal\\%%4
\hline
\multicolumn{1}{c|}{$n\equiv 3,5,7,9,11,13 \mod 14$} & $7n$ &$8(n-1)$&
minimal\\%%4
\hline
\multicolumn{1}{c|}{$n\equiv 1 \mod 14, n\neq 15$} &$n$ & $8(n-1)/7$&  7 times
the minimal \\%%5
\hline
\multicolumn{1}{c|}{$n= 15$} &15 & 16& minimal \\%%5
\hline
\end{tabular}}\label{42summary}
\end{table}
\end{center}
Table \ref{42summary} summarizes the results given in this section. It shows
when the necessary
conditions are sufficient for (2,4) GDDs with minimal indices. Further, the
table indicates when the results show the necessary conditions are
sufficient for configuration $(2,4)$ GDDs with near minimal, seven times the
minimal possible or fourteen times the minimal possible indices. 

%%%%%%%%%%%%%%%%%%%%%%%%%%%%%%%%%%%
%%%%%%%%%%%%%%%%%%%%%%%%%%%%%5%%%%%
%CONFIG 5,1
%%%%%%%%%%%%%%%%%%%%%%%%%%%%%%%
%%%
\section{GDDs with Configuration (1,5)}

In this section we focus on the minimal indices for configuration $(1,5)$
GDD$(n,$$2,$ $6;$$\lo,\lt)$. 
Hurd and Sarvate gave a construction for configuration $(1,k)$
GDD$(n,2,k+1;\lambda_1,\lambda_2)$ using a BIBD($n,k,\Lambda$)s \cite{5gdd14}.
We repeat their result here:

\begin{theorem}[\cite{5gdd14}]\label{kplus1gdd}
 The existence of a BIBD$(n,k,\Lambda)$ implies the existence of a
configuration $(1,k)$ GDD$(n,2,k+1;\lambda_1,\lambda_2)$ with
$\lambda_1=\Lambda n$ and $\lambda_2=2\Lambda(n-1)/(k-1)$.
\end{theorem}

\begin{center}
\begin{table}
\caption{Existence of BIBD$(n,5,\la)$ and Resulting Configuration $(1,5)$ GDDs.}
\begin{tabular}{l|l|l}

 \hline
 \multicolumn{1}{c|}{BIBD}  &Existence  & Resulting GDD  \\
\hline\hline
 \multicolumn{1}{c|}{$(n,5,1)$}  &$n\equiv 1,5\mod 20$  & $GDD(n,2,6;n,(n-1)/2)$
 \\
\hline
 \multicolumn{1}{c|}{$(n,5,2)$}  &$n\equiv 1,5\mod 10, n\neq 15$  &
$GDD(n,2,6;2n,n-1)$ 
\\
\hline
 \multicolumn{1}{c|}{$(n,5,4)$}  &$n\equiv 0,1\mod 10, n\neq 10,160,190$  &
$GDD(n,2,6;4n,2(n-1))$
 \\
 \hline
 \multicolumn{1}{c|}{$(n,5,5)$}  &$n\equiv1\mod 4$  & $GDD(n,2,6;5n,5/(2(n-1)))$
 \\
 \hline
 \multicolumn{1}{c|}{$(n,5,10)$}  &$n\equiv1\mod 2$  & $GDD(n,2,6;10n,5(n-1))$ 
\\
 \hline
 \multicolumn{1}{c|}{$(n,5,20)$}  & All $n$  & $GDD(n,2,6;20n,10(n-1))$  \\
 \hline
\end{tabular}\label{BIBDGDDs}
\end{table}
\end{center}

Further, in \cite{BIBD5} Hanani showed the existence of some classes of
BIBD$(n,5,\lambda)$.  Using his result and Theorem \ref{kplus1gdd} we obtain the
following $(1,5)$ configuration GDD$(n,2,6;$$\lo,$$\lt)$s summarized in Table 
\ref{BIBDGDDs}. 

However, this construction does not always give optimal values of $\lambda_1$
and
$\lambda_2$. By Theorem \ref{betahurd}, we have the following relation between
\lo\; and \lt. 

\begin{corollary}\label{fiveonelambdas}
 For any configuration $(1,5)$ GDD$(n,2,6;$$\lambda_1,$$\lambda_2)$ we have\\
$\lambda_2=\frac{\lambda_1(n-1)}{2n}$. 
\end{corollary}
From Theorem \ref{fiveonelambdas} we see that for some value of $n$ the minimum
possible indices are $(\lo,\lt)$. As in the other two configurations, we may
find the minimal indices by Theorem \ref{fiveonelambdas} and Theorem \ref{thm1}.
Further, any other GDD with configuration $(1,5)$ will have indices
$(w\lo,w\lt)$ for some positive integer $w$. The minimal indices are summarized
in the next theorem. 
\begin{theorem}\label{lambdaTable15}
The minimal indices $(\lo,\lt)$ for any configuration $(1,5)$
GDD$(n,2,6;$ $\lo,\lt)$ summarized in Table \ref{15table}.

\begin{table} 
\caption{Summary of Minimal Indices for Configuration $(1,5)$}
\begin{center}
\begin{tabular}{l|r|r}
 \hline%%0
 \multicolumn{1}{c|}{$n$} &$\lambda_1$ &$\lambda_2$  \\%%1
\hline
\multicolumn{1}{c|}{$n\equiv 0,6,10,11,15,16\mod{20}$} & $2n$  &$(n-1)$\\%%1
 \hline
\multicolumn{1}{c|}{$n\equiv 1,5\mod{20}$} & $n$ &$(n-1)/2$\\%%1
 \hline
\multicolumn{1}{c|}{$n\equiv 2,4,8,12,14,18\mod{20}$} & $10n$ &$5(n-1)$\\%%1
 \hline
\multicolumn{1}{c|}{$n\equiv 3,7,9,13,17,19\mod{20}$} & $5n$ &$5(n-1)/2$\\%%1
 \hline
\end{tabular}\label{15table}
\end{center}
\end{table}

\end{theorem}
\begin{proof}

By Theorem \ref{fiveonelambdas}, we have that
$\lambda_2=\displaystyle\frac{\lambda_1(n-1)}{2n}$. This implies that if
$n\equiv 1\mod 2$ then $\lo$ must be a multiple of $n$ and $\lt$ must be a a
multiple of $(n-1)/2$. However, if $n\equiv 11,15\mod 20$ then the indices given
do not give an even number of blocks which is required by Theorem
\ref{evenblocks}. So for $n\equiv 11,15\mod 20$, if we take two times the
minimum possible indices, the number of blocks will be integer valued implying
$(\lo,\lt)=(2n,(n-1))$. Also, using the given indices for $n\equiv 3,7,9\mod 10$
results in a non-integer value for the number of blocks given by Theorem
\ref{thm1}. Thus we must take 5 times these, so the minimal indices are
$(\lo,\lt)=(5n,5(n-1)/2)$. Finally, if $n\equiv 1,5\mod 20$, the necessary
conditions in Theorem \ref{thm1} are met.

If $n\equiv 0\mod 2$, Theorem \ref{fiveonelambdas} tells us that $\lo$ must be a
multiple of $2n$ and $\lt$ must be a multiple of $n-1$. However if $n\equiv
2,4,8\mod 10$, then these values give a non-integer value for the number of
blocks. If we take 5 times these indices then the necessary condition in Theorem
\ref{thm1} is satisfied, and so the minimal indices are
$(\lo,\lt)=(10n,5(n-1))$. Notice that for $n\equiv 0,6\mod 10$, the given
indices are $(\lo,\lt)=(2n,n-1)$ which are the minimum possible. 
\end{proof}

%%%%%%%%%%%%%%%%%%%%%%%%%%%%%%%%%%%
%%%%%%%%%%%%%%%%%%%%%%%%%%%%%5%%%%%
%CONFIG 5 1 
%%%%%%%%%%%%%%%%%%%%%%%%%%%%%%%%%%
\section{Constructing Configuration (1,5) GDDs}
In this section we focus on constructing $(1,5)$ GDDs with minimal indices.
Theorem \ref{kplus1gdd} gives us the following results. 
\begin{corollary}
 The necessary conditions are sufficient for the existence of a configuration
$(1,5)$ GDD$(n,2,6;n,(n-1)/2)$ for $n\equiv 1,5\mod 20$.
\end{corollary}

\begin{corollary}
 The necessary conditions are sufficient for the existence of a configuration
$(1,5)$ GDD$(n,2,6;2n,n-1)$ for $n\equiv 11,15\mod 20, n\neq 15.$
\end{corollary}
Notice that in the previous two constructions, the design is minimal. 
We use a resolvable BIBD$(n,5,4)$ in the
following construction. In \cite{CRCres}, it is given that a resolvable
BIBD$(n,5,4)$ exists for $n\equiv 0 \mod 10$ except for $n= 10,160,190$.

%%%%%%%%%%%%% 0 mod 10
\begin{theorem}
 Let $n\equiv 0 \mod{10}, n\neq 10,160,190$. Then the necessary conditions are
sufficient 
for the existence of a configuration $(1,5)$
GDD$(n,2,6;$ $2n, n-1)$. 
\end{theorem}
\begin{proof}
  Let $n\equiv 0 \mod{10}, n\neq 10,160,190$. Assume the two groups are
$A=\{1,2,\dots,n\}$ and $A'=\{1',2',\dots,n'\}$. There exists a 
RBIBD($n,5,4$)
with
$b=n(n-1)/5$ blocks, and each point appearing $r=(n-1)$ times. Let $D$ be such
a design on $A$ with parallel classes $C_1,C_2,\dots, C_{n-1}$. Construct a
graph $G$ in the following manner. For $j=1,2,\dots,n-1$, create the bipartite
graph $G_j$ with bipartitions $G_j^1$ and $G_j^2$ where $V(G_j^1)$ are the
blocks of $C_j$ and $V(G_j^2)$ are the points in $A'$. Each of the first $n/10$
vertices in $G_j^1$ are adjacent to the vertices in $G_j^2$ that correspond to
the first $n/10$ blocks of $C_j$. Each of the last $n/10$ vertices in $G_j^2$
are adjacent to the vertices in $G_j^2$ that correspond to the last $n/10$
blocks of $C_j$. Thus each vertex in $G_j^1$ has degree $n/2$ and each vertex
in $G_j^2$ has degree $n/10$. This creates half of the desired blocks in the
GDD. To obtain the other half, let $D$ be an RBIBD$(n,5,4)$ on $A'$ and repeat
the construction. This time, each of the first $n/10$ vertices in $G_j'$ will be
adjacent to the vertices in $G_j^2$ that correspond to the last $n/10$ blocks
of $C_j$, and each of the last $n/10$ vertices of $G_j^1$ will be adjacent to
the vertices in $G_j^2$ that correspond to the first $n/10$ blocks of $C_j$. 

In the design, each pair appears four times and will be matched $n/2$ times.
Now consider a second pair of associates $\{x,y'\}$ where $x\in A$ and $y'\in
A'$. The points $x$ and $y$ appear in every parallel class exactly once. So
for each $G_j$, if $x$ and $y$ are both in the same half of $A$ (either in the
first $n/10$ blocks of $C_j$ or the last $n/10$ blocks of $C_j$) then
$\{x,y'\}$ appears once in the first part of the construction and zero times in
the second part. If $x$ and $y$ were in different halves of $A$, then
$\{x,y'\}$ appears once in the second part of the construction and zero times
in the first part. Therefore $\{x,y'\}$ appears exactly once per $G_j$. 
Thus $\lambda_2$ is
the number of parallel classes or $n-1$. 
\end{proof}
A near parallel class is a partial parallel class missing a single point. A
near - resolvable design NRB($n,k,k-1$) is a BIBD($n,k,k-1$) with the
property that the blocks can be partitioned into near parallel classes. For
such a design, every point is absent from exactly one class. The necessary
condition for the existence of an NRB($v,k,k-1$) is $v\equiv 1 \mod k$. It is
known that the necessary condition is sufficient for the existence of a
NRB($v,k,k-1$) if $k\leq 7$ (see \cite{CRCres}). We use near resolvable designs
in the
following
construction. 
%%%%%%%%%%%%% 6 mod 10
\begin{theorem}
 Let $n\equiv 6 \mod{10}$. Then the necessary conditions are
sufficient 
for the existence of a configuration $(1,5)$
GDD$(n,2,6;2n, n-1)$. 
\end{theorem}
\begin{proof}
  Let $n\equiv 6 \mod{10}$, and the two groups have point sets $A=\{1,2,\dots,
$ $n\}$ and $A'=\{1',2',\dots, n'\}$. Since $n\equiv 6 \mod{10}$, there exists a
NRB($n,5,4$). It has $n$
near parallel classes with $(n-1)/5$ blocks in them each. Let $D$ be such a
design on the point set of $A$, and resolve the blocks of $D$ into near
parallel classes $C_1,C_2,\dots C_n$ where $C_i$ misses point $i$. Construct a
graph $G$ in the following manner. For $j=1,2,\dots n/2$, create the complete
bipartite graph $G_j$ with bipartitions $G_j^1$ and $G_j^2$ where $V(G_j^1)$
are the blocks of $C_j$ and $V(G_j^2)$ are the points $\{1',2',\dots, n/2'\}$.
For $j=n/2+1,\dots,n$, create the complete bipartite graph $G_j$ with
bipartitions $G_j^1$ and $G_j^2$ where $V(G_j^2)$ are the points
$\{(n/2+1)',(n/2+2)',\dots, n'\}$. This creates half of the desired blocks. To
get the other half, let $D$ be the NRB$(n,5,4)$ on $A'$ and repeat the
construction with $V(G_j^2)$ being the points $\{1,2,\dots,n/2\}$ for
$j=n/2+1,\dots,n$, and $V(G_j^2)$ being the points $\{(n/2)+1,\dots,n\}$ for
$j=1,2,\dots, n/2$. 

Consider a pair of first associates. It will appear $4(n/2)=2n$ times in a
block of size 6. Now consider a pair of second associates where
$x\in A$ and $y'\in A'$. If $x\in \{1,2,\dots, n/2\}$ and $y'\in \{1', 2',
\dots,
(n/2)'\}$ then $\{x,y'\}$ will appear $(n/2)-1$ times in the first part of the
construction
and $n/2$ times in the second. It is the same case if $x\in
\{n/2+1,n/2+2,\dots, n\}$ and $y'\in \{(n/2+1)',(n/2+2)',\dots, n'\}$. If $x\in
\{1,2,\dots, n/2\}$ and $y'\in \{(n/2+1)',(n/2+2)',\dots, n'\}$, then $\{x,y'\}$
will
appear 
$n/2$ times in the first part and $n/2-1$ times in the second part. It is the
same case if  $x\in \{n/2+1,n/2+2,\dots, n\}$ and  $y'\in \{1', 2', \dots,
n'\}$.
Thus $\lt=n-1$. 
\end{proof}
Note that we have constructed minimal GDDs for $n\equiv 0,1,5,6\mod 10$ (for all
but a few values). Recall that a near-minimal design is one that has exactly
twice the minimal indices. By Theorem \ref{kplus1gdd}, the necessary conditions
are sufficient for the existence of a near minimal GDD$(n,2,6;\lo,\lt)$ for
$n\equiv 2,3,4,7,8,9\mod 10$. We may construct a minimal GDD$(n,2,6;\lo,\lt)$
for $n\equiv 3,7,9\mod 10$ given the existence of a 5-resolvable
BIBD$(n,5,10)$. 

\begin{theorem}
The existence of a $5$-resolvable BIBD$(n,5,10)$ implies the existence of a
configuration $(1,5)$ GDD$(n,2,6;5n,5(n-1)/2)$ for $n\equiv 3,7,9 \mod{10}$. 
\end{theorem}
\begin{proof}
  Let $n\equiv 3,7,9 \mod{10}$ and assume there exists a $5$-resolvable
BIBD $(n,5,10)$. Assume the two groups are $A=\{1,2,3,\dots, n\}$ and
$A'=\{1',2',3',$ $\dots, n'\}$ and let $D$ be such a design on point set $A$.
Resolve the blocks of $D$ into 5-parallel classes $C_1,C_2,\dots, C_{n-1/2}$,
each having $n$ blocks. Construct a graph $G$ in the following manner. For
$j=1,2,\dots, (n-1)/4$, create the complete bipartite graph $G_j$ with
bipartitions $G_j^1$ and $G_j^2$ where $V(G_j^1)$ are the blocks of $C_j$ and
$V(G_j^2)$ are the odd numbers in $A'$. For $j=(n-1)/4+1,\dots, (n-1)/2$,
create the complete bipartite graph $G_j$ with bipartitions $G_j^1$ and $G_j^2$
where $V(G_j^1)$ are the blocks of $C_j$ and $V(G_j^2)$ are the even numbers in
$A'$. This creates half of the desired blocks. To get the other half, let $D$
be a 5-RBIBD$(n,5,10)$ on $A'$ and repeat the construction with $V(G_j^2)$
being the even numbers in $A$ for $j=1,2,\dots, (n-1)/4$ and $V(G_j^2)$ being
the odd numbers in $A$ for $j=(n-1)/4+1,\dots,(n-1)/2$. 

Consider a pair of first associates. It will appear 10 times in $D$. Therefore,
in the given construction it will appear $5n$ times in a block of size 6. Now
consider a pair of second associates $\{x,y'\}$. In each part of the
construction, this pair appears $5(n-1)/4$ times, thus it appears a total of
$5(n-1)/2$ times. 
\end{proof}

\subsection{Summary of Minimality}
\begin{center}
\begin{table}[here]
\caption{Summary of Constructions and Minimality for Configuration $(1,5)$}
\scalebox{0.8}{
\begin{tabular}{l|r|r|| c}

 \hline%%0
 \multicolumn{1}{c|}{$n$} &$\lambda_1$ &$\lambda_2$  \\%%1
\hline
\multicolumn{1}{c|}{$n\equiv 0,10,11,15,6,16\mod{20}$, $n\neq 10,15,160,190$} &
$2n$ &$(n-1)$& minimal\\%%1
 \hline
\multicolumn{1}{c|}{$n\equiv 1,5\mod{20}$} & $n$ &$(n-1)/2$& minimal\\%%1
 \hline
\multicolumn{1}{c|}{$n\equiv 2,4,8,12,14,18\mod{20}$} & $10n$ &$5(n-1)$&
near-minimal\\%%1
 \hline
\multicolumn{1}{c|}{$n\equiv 3,7,9,13,17,19\mod{20}$} & $5n$ &$5(n-1)/2$&
near-minimal\\%%1
 \hline
\end{tabular}\label{construction15sum}
}
\end{table}
\end{center}

We conclude this section with a summary of the GDDs we have constructed, and
their minimality found in Table \ref{construction15sum}.

\bibliographystyle{amsplain}
\bibliography{thesis.bib}

\end{document}